\documentclass[12pt]{amsart} 
\usepackage{graphpap}  
\begin{document}  

\newcommand{\ii}{\hookrightarrow}  
\newcommand{\hh}{\ensuremath{\mathbb H}}  
\newcommand{\ee}{\ensuremath{\mathbb E}}  
\newcommand{\rr}{\ensuremath{\mathbb R}}
\newcommand{\qq}{\ensuremath{\mathbb Q}}
\newcommand{\zz}{\ensuremath{\mathbb Z}}
\newcommand{\cc}{\ensuremath{\mathcal{C}}}
\newcommand{\ff}{\ensuremath{\mathcal{F}}}
\newcommand{\uu}{\ensuremath{\mathcal{U}}}

\newcommand{\al}{\ensuremath{\alpha}}
\newcommand{\bb}{\ensuremath{\beta}}
\newcommand{\be}{\ensuremath{\beta}}
\newcommand{\g}{\ensuremath{\gamma}}
\newcommand{\G}{\ensuremath{\Gamma}}
\newcommand{\de}{\ensuremath{\delta}}
\newcommand{\D}{\ensuremath{\Delta}}
\newcommand{\ep}{\ensuremath{\varepsilon}}
\newcommand{\f}{\ensuremath{\varphi}}
\newcommand{\s}{\ensuremath{\sigma}}  
\newcommand{\la}{\ensuremath{\lambda}}
\newcommand{\w}{\ensuremath{\omega}}  
\newcommand{\W}{\ensuremath{\Omega}}

\newcommand{\p}{\ensuremath{\partial}}  
\newcommand{\M}{\ensuremath{\widetilde{M}}}  
\newcommand{\e}{\ensuremath{\epsilon}}  
\newcommand{\pp}{\ensuremath{\parallel}} 
\newcommand{\inv}{\ensuremath{^{-1}}} 
\newcommand{\se}{\ensuremath{\subseteq}}
\newcommand{\lra}{\ensuremath{\longrightarrow}}
\newcommand{\lla}{\ensuremath{\longleftarrow}}
\newcommand{\pf}{{\noindent\bf Proof.\ }}
\newcommand{\noi}{\noindent}
\newcommand{\ds}{\displaystyle}

\newcommand{\go}{\ensuremath{\G^{(0)}}}
\newcommand{\fff}{\ensuremath{\bar{f}}}
\newcommand{\?}{{\bf ?????}}
\newcommand{\bu}{$\bullet$\ }

\newtheorem{ttt}{Theorem}  
\newtheorem{lll}[ttt]{Lemma}
\newtheorem{slll}[ttt]{Sublemma}
\newtheorem{ddd}[ttt]{Definition}  
\newtheorem{ppp}[ttt]{Proposition}  
\newtheorem{ccc}[ttt]{Corollary}
\newtheorem{nnn}[ttt]{Not a theorem}

\newtheorem{ttt-in-text}{Theorem}
%\renewcommand{\baselinestretch}{1.3} 
%\small\normalsize  

\title[The  Baum-Connes conjecture  for hyperbolic groups]
{The  Baum-Connes conjecture for hyperbolic groups}
\maketitle
\centerline{Igor Mineyev and Guoliang 
Yu\footnote{ The second author is
 partially supported by
NSF and MSRI. }} 
\date{}
\begin{abstract}
{We prove the Baum-Connes conjecture for hyperbolic groups and \\ their
subgroups. } 
\end{abstract}
%\maketitle 
%\renewcommand{\thefootnote}{\fnsymbol{footnote}}
%\footnote{\textit{Date}: The second author is
%partially supported by
%NSF and MSRI }  
%\footnote[0]{1991 \textit{Mathematics Subject Classification}:
%}
%\setcounter{section}{-1} 

\section{Introduction.}

The Baum-Connes conjecture states that, for a discrete group $G$, 
the K-homology groups
of the classifying space for proper $G$-action
is isomorphic to  the K-groups of the reduced
group $C^{\ast}$-algebra of  $G$~\cite{BC, BCH}. 
A positive answer to the Baum-Connes conjecture would provide a complete
solution to the problem of computing higher indices of elliptic 
operators on compact manifolds.
The rational injectivity part of the Baum-Connes conjecture implies 
the Novikov conjecture on homotopy invariance of higher signatures. 
 The Baum-Connes conjecture  also implies the Kadison-Kaplansky 
conjecture that for $G$ torsion free there exists no non-trivial
projection in the reduced group $C^{\ast}$-algebra associated to $G$.
In~\cite{HK}, Higson and Kasparov prove the Baum-Connes conjecture for groups
acting properly and isometrically on a Hilbert space.
In a recent remarkable work, Vincent Lafforgue proves the Baum-Connes
conjecture for strongly bolic groups with property RD~\cite{L1, L2, L3}. 
In particular, this implies the Baum-Connes conjecture for
the fundamental groups of strictly
negatively curved compact manifolds. In~\cite{CM}, Connes and Moscovici
prove the rational injectivity
part of the Baum-Connes conjecture for hyperbolic groups using cyclic
cohomology method. In this paper, we exploit Lafforgue's work
to prove the Baum-Connes conjecture for hyperbolic groups and their
subgroups.

The main step in the proof is the following theorem.
\setcounter{ttt-in-text}{17}
\addtocounter{ttt-in-text}{-1}
\begin{ttt-in-text}
Every hyperbolic group $G$ admits a metric $\hat{d}$ with the following properties.
\begin{itemize}
\item [(1)] $\hat{d}$ is $G$-invariant, i.e.
$\hat{d}(g\cdot x,g\cdot y)=\hat{d}(x,y)$ for all $x,y,g\in G$.
\item [(2)] $\hat{d}$ is quasiisometric to the word metric.
\item [(3)] The metric space $(G,\hat{d})$ is weakly geodesic and strongly bolic.
\end{itemize}
\end{ttt-in-text}

This paper is organized as follows. In section 2, we recall the concepts of
hyperbolic groups and bicombings. In section 3, we introduce a distance-like 
 function $r$ on a hyperbolic group and study its  basic properties.
In section 4, we prove that $r$ satisfies certain distance-like 
inequalities. In section 5, we construct a metric $\hat{d}$ on a hyperbolic
group and prove
Theorem~\ref{metric} stated above.
In section 6, we combine Lafforgue's work and Theorem 17 to
 prove the Baum-Connes conjecture for hyperbolic groups
and their subgroups.  

After this work was done, we learned from 
 Vincent Lafforgue that he has
independently proved the Baum-Connes conjecture for hyperbolic groups
by a different and elegant method~\cite{L4},
and we also learned from Michael Puschnigg that he has independently
proved the Kadison-Kaplansky conjecture for hyperbolic groups
using  a beautiful  local cyclic homology method~\cite{p}.
It is our pleasure to thank both of them for bringing their work to our
attention.  

\section{Hyperbolic groups and bicombings.}

In this section, we recall the concepts of hyperbolic groups and 
bicombings.

\subsection{Hyperbolic groups.}
\label{s_hyp}

Let $G$ be a finitely generated group.
 Let $S$ be a  finite generating 
set for $G$. Recall that the Cayley graph  of $G$ with respect to $S$
is the graph $\G$  satisfying the following conditions:  
\begin{itemize}
\item [(1)] the set of vertices in $\G$, denoted by $\go$, is $G$;
\item [(2)] the set of edges is $G\times S$, where each edge
  $(g,s)\in G\times S$ spans the vertices $g$ and $gs.$
\end{itemize}
  
We endow  $\G$  with the {\sf path metric}~$d$
induced by assigning length 1 to each
edge.
Notice that $G$ acts freely, isometrically  and cocompactly on $\G$. 
 A {\sf geodesic path} in $\G$ is a shortest edge path. 
The restriction of the path metric $d$ to $G$ is called the {\sf word metric}.

A finitely generated group $G$ is called {\sf hyperbolic}\,
if
there exists a constant $\de\ge 0$ such that all the geodesic triangles
in $\G$ are {\sf $\de$-fine} in the following sense:
if $a$, $b$, and $c$ are vertices in $\G$,
$[a,b]$, $[b,c]$, and $[c,a]$ are geodesics from $a$ to $b$,
from $b$ to $c$, and from $c$ to $a$, respectively, and
points $\bar{a}\in [b,c]$, $v,\bar{c}\in [a,b]$,
$w,\bar{b}\in [a,c]$ satisfy
$$d(b,\bar{c})=d(b,\bar{a}),\quad d(c,\bar{a})=d(c,\bar{b}),\quad
d(a,v)=d(a,w)\le d(a,\bar{c})=d(a,\bar{b}),$$
then $d(v,w)\le \de$.

The above definition of hyperbolicity  does not depend on the choice of
the  finite generating set $S$.
 See~\cite{G, ABC}
 for other equivalent definitions.

For vertices $a$, $b$, and $c$ in $\G$,
the {\sf Gromov product} is defined by
$$(b|c)_a := d(a,\bar{b})=  d(a,\bar{c})=
  \frac{1}{2}\Big[d(a,b)+d(a,c)-d(b,c)\Big].$$
The Gromov product can be used to 
 measure the degree of cancellation in the multiplication 
of group elements in $G$. 

\subsection{Bicombings.}
Let $G$ be a finitely generated group. Let $\G$ be a
Cayley graph with respect to a finite generating set. 
A {\sf bicombing} $p$ in $\G$ is a function assigning
to each ordered pair $(a,b)$ of vertices in $\G$
an oriented edge-path $p[a,b]$ from $a$ to $b$.
A bicombing $p$ is called {\sf geodesic}\, if
each path $p[a,b]$ is geodesic, i.e. a shortest edge path.
A bicombing $p$ is {\sf $G$-equivariant} if
$p[g\cdot a, g\cdot b]= g\cdot p[a,b]$
for each $a,b\in\go$ and each $g\in G$.

\section{Definition and properties of $r(a,b)$.}

The purpose of this section is to introduce a distance-like function $r$
on a hyperbolic group and study its basic properties.

Let $G$ be a hyperbolic group and $\G$ be a Cayley graph of $G$ with respect
to a finite  generating set.
 We endow $\G$ with  the  path metric $d$, 
 and identify $G$ with $\go$,  the set of vertices of $\G$.
Let $\de\ge 1$ be a positive integer such that 
all the geodesic triangles in $\G$
are $\de$-fine.

The {\sf ball} $B(x,R)$ is the set of all vertices at distance
at most $R$ from the vertex $x$. The {\sf sphere} $S(x,R)$ is the set
of all vertices at distance $R$ from the vertex $x$.
Pick an equivariant geodesic bicombing~$p$ in~$\G$. By $p[a,b](t)$ we denote
the point on the geodesic path $p[a,b]$ at distance $t$ from~$a$.
Recall that $C_0 (G,\qq)$ is the space of all
 $0$-chains (in $G=\go$)  with coefficients in
$\qq$.
Endow $C_0(G,\qq)$ with the $\ell^1$-norm $|\cdot|_1$.
We identify $G$ with the standard basis of $C_0(G,\qq)$.
Therefore the left action of $G$ on itself induces a left action
on~$C_0(G,\qq)$.

First we recall several  constructions from~\cite{M3}.

For $v, w\in G$, the {\sf flower} at $w$ with respect to $v$ is defined to be
$$ Fl(v, w) := S(v, d(v,w))\cap B(w,\de)\subseteq G.$$

For each $a\in G$, we define $pr_a: G\rightarrow G$ by:
\begin{itemize}
\item [(1)] $pr_a (a) := a;$
\item [(2)] if $b\neq a$, $pr_a(b) := p[a,b](t)$, where $t$ is the largest
integral multiple of $10\de$ which  is strictly less than $d(a,b)$.
\end{itemize}
 
Now for each pair $a,b\in G$, we define a $0$-chain $f(a,b)$ in 
$G$ inductively on the distance $d(a,b)$ as follows:  
\begin{itemize}
\item [(1)] if $d(a,b)\leq 10\de$, $f(a,b) :=b$;
\item [(2)] if $d(a,b) >10\de$ and $d(a, b)$ is not an integral
multiple of $10\de$, let $f(a,b) := f(a, pr_a (b));$  
\item [(3)] if $d(a,b)> 10\de$ and $d(a,b)$ is an integral
 multiple of $10\de$, 
 let
$$ f(a,b) := \frac{1}{\# Fl(a,b)}\sum_{x\in Fl(a,b)} f(a, pr_a (x)).$$
\end{itemize}
\begin{ppp}[\cite{M3}]
\label{c-c-c}
The function $f:G\times G\to C_0(G,\qq)$
defined above satisfies the following conditions.
\begin{itemize}
\item [(1)] For each $a,b\in G$, $f(b,a)$ is a {\sf convex combination},
i.e. its coefficients are non-negative and sum up to~1.
\item [(2)] If $d(a,b)\ge 10\de$, then
$supp\,f(b,a)\subseteq B(p[b,a](10\de),\de)\cap S(b,10\de)$.
\item [(3)] If $d(a,b)\le 10\de$, then $f(b,a)=a$.
\item [(4)] $f$ is $G$-equivariant, i.e.
$f(g\cdot b, g\cdot a)= g\cdot f(b,a)$
for any $g,a,b\in G$.
\item [(5)] There exist constants $L\ge 0$ and $0\le\la<1$
such that, for all $a,a',b\in G$,
$$\Big|f(b,a)-f(b,a')\Big|_1\le L\,\la^{(a|a')_b}.$$
\end{itemize}
\end{ppp}
\setlength{\unitlength}{.7cm} 
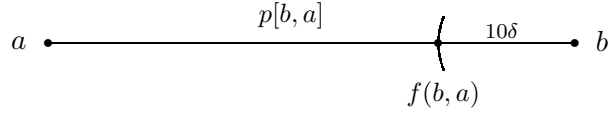
\begin{figure}[h]
  \begin{center}
   \begin{picture}(12,5.2)
\put(1,2.5){\circle*{0.15}}
\put(.3,2.4){\small$a$}

\put(11,2.5){\circle*{0.15}}
\put(11.4,2.35){\small$b$}

\put(5,2.9){\footnotesize$p[b,a]$}
\put(9.3,2.6){\tiny $10\de$}

%from a to b
\qbezier(1,2.5)(6,2.5)(11,2.5)

%f(a,b)
\qbezier(8.5,2)(8.3,2.5)(8.5,3)
\put(7.8,1.4){\footnotesize$f(b,a)$}
\put(8.41,2.5){\circle*{0.15}}
   \end{picture}
  \end{center}
 \caption{Convex combination $f(b,a)$.}
\end{figure}
Let $\w_7$ be the number of elements in a ball of radius $7\de$ in $G$.
For each $a\in G$, a 0-chain $star(a)$ is defined by
$$star(a):=\frac{1}{\w_7}\sum_{x\in B(a,7\de)} x.$$
This extends to a linear operator $star: C_0(G,\qq)\to C_0(G,\qq)$.
Define the 0-chain $\fff(b,a)$ by $\fff(b,a):=star\big( f(b,a) \big)$.

The main reason for introducing $\fff$ is that $\fff$ has better
cancellation properties than $f$ (compare Proposition~\ref{c-c-c}(5) with
Proposition~\ref{la'}(5) and~\ref{la'}(6) below).
These cancellation properties play key roles in this paper.

\begin{ppp}[\cite{M3}]
\label{la'}
The function $\fff:G\times G\to C_{0}(G,\qq)$ defined above
satisfies the following conditions.
\begin{itemize}
\item [(1)] For each $a,b\in G$, $\fff(b,a)$ is a convex combination.
\item [(2)] If $d(a,b)\ge 10\de$, then
$supp\,\fff(b,a)\subseteq B(p[b,a](10\de),8\de)$.
\item [(3)] If $d(a,b)\leq 10\de$, then $supp\, \fff(b,a) \subseteq 
B(a, 7\de).$ 
\item [(4)] $\fff$ is $G$-equivariant, i.e.
$\fff(g\cdot b, g\cdot a)= g\cdot \fff(b,a)$
for any $g,a,b\in G$.
\item [(5)] There exist constants $L\ge 0$ and $0\le\la<1$
such that, for all $a,a',b\in G$,
$$\Big|\fff(b,a)-\fff(b,a')\Big|_1\le L\,\la^{(a|a')_b}.$$
\item [(6)] There exists a constant $0\le\la'<1$ such that if $a,b,b'\in G$
satisfy $(a|b)_{b'}\le 10\de$ and $(a|b')_b\le 10\de$, then
$\Big|\fff(b,a)-\fff(b',a)\Big|_1\le 2\la'$.
\item [(7)] Let $a,b,c\in G$, $\gamma$ be a geodesic path from $a$ to $b$,
and let
$$c\in N_G(\gamma, 9\delta):=\{x\in G\ \big|\  d(x,\gamma)\leq 9\delta\}.$$
Then $supp (\bar{f}(c,a))\subseteq N_G (\gamma, 9\delta).$ 
\end{itemize}
\end{ppp}

\begin{ddd}
\label{def-r}
For each pair of vertices $a,b\in G$, a rational number $r(a,b)\ge 0$ is
defined inductively on $d(a,b)$ as follows.
\begin{itemize}
\item $r(a,a):=0$.
\item If $0< d(a,b)\le 10\de$,
let $r(a,b):= 1$.
\item If $d(a,b)> 10\de$, let
$r(a,b):=r\big(a,\fff(b,a)\big)+ 1$, where
$r\big(a,\fff(b,a)\big)$ is defined by linearity in the second variable.
\end{itemize}
\end{ddd}

The function $r$ is well defined by Proposition~\ref{la'}(2). Also, $r(a,b)$
is well defined when $b$ is a 0-chain, by linearity.

Let $\qq_{\ge 0}$ denote the set of all non-negative rational numbers.

\begin{ppp}
\label{N}
For the function $r: G\times G\to \qq_{\ge 0}$ defined above,
there exists $N\ge 0$ such that,
for all $a,b,b'\in G$,
$$\big|r(a,b)-r(a,b')\big|\le d(b,b')+N.$$
\end{ppp}
\pf 
Up to the $G$-action, there are only finitely many triples of
vertices $a$, $b$, $b'$, satisfying $d(a,b)+d(a,b')\le 40\de$,
hence there exists a uniform bound $N'$ for the norms
$$\big|r(a,b)-r(a,b')\big|$$
for such vertices $a$, $b$, $b'$.
Let $\la'$ be the constant from Proposition~\ref{la'}(6) and
pick $N$ large enough so that
\begin{equation}
\label{N-def}
N'\le N\qquad\mbox{and}\qquad \la'\cdot[27\de+N]\le N.
\end{equation}
We shall prove the inequality in Proposition 4 by induction on $d(a,b)+d(a,b')$.

If $d(a,b)+d(a,b')\le 40\de$, then
$$\big|r(a,b)-r(a,b')\big|\le N'\le N\le d(b,b')+N$$
just by the choices of $N'$ and $N$.
We assume now that $d(a,b)+d(a,b')> 40\de$. Consider the following two
cases.

\noi{\underline{\bf Case~1.}
{\it \ $(a|b')_b> 10\de$\quad or\quad $(a|b)_{b'}> 10\de$.}

\setlength{\unitlength}{1 cm} 
\begin{figure}[h]
  \begin{center}
   \begin{picture}(8,4.5)
%\graphpaper[1](0,0)(8,4.5)

\put(0.5,1.5){\circle*{0.1}}
\put(0.1,1.5){\small$a$}

\put(7,4){\circle*{0.1}}
\put(7.2,4){\small$b$}

\put(6,0){\circle*{0.1}}
\put(6.2,-0.1){\small$b'$}

%from a to b
\qbezier(0.5,1.5)(4.5,2)(7,4)
\put(2.8,2.3){\footnotesize$p[b,a]$}
%from a to b'
\qbezier(0.5,1.5)(4.5,2)(6,0)
%from b to b'
\qbezier(7,4)(4.5,2)(6,0)
\put(6,1.1){\small$\g$}
%inscribed triangle
\put(4.42,2.49){\circle*{0.1}}
\put(5.455,1.7){\circle*{0.1}}
\put(4.47,1.18){\circle*{0.1}}

\put(6.2,2.8){\tiny$10\de$}

\put(5.5,1.97){\circle*{0.1}}
\put(5.7,1.9){\small$v'$}

\put(4.7,2.6){\circle*{0.1}}
\put(5,2.48){\small$v$}

\put(4.7,3.9){\circle*{0.1}}
\put(4.2,3.9){\small$x$}

%from v to x
\qbezier(4.7,2.6)(4.7,3.25)(4.7,3.9)
%from v to v'
\qbezier(4.68,2.6)(5.09,2.285)(5.5,1.97)
   \end{picture}
  \end{center}
 \caption{\label{case1} Proposition~\ref{N}, Case~1.}
\end{figure}
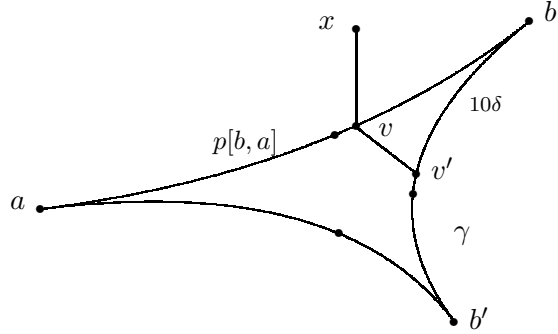

Assume, for example, that $(a|b')_b> 10\de$. 
Then
$d(a,b)> 10\de$, hence, by definition,
$$r(a,b)=r\big(a,\fff(b,a)\big)+ 1.$$
By Proposition 2(2), we have
 $supp\,\fff(b,a)\se B(v,8\de)$, where $v:=p[b,a](10\de)$.
Also, $(a|b')_b>10\de$ implies  $d(b,b')> 10\de$.
Hence  there exists a geodesic $\g$
between $b$ and $b'$, and a vertex $v'$ on $\g$ with
$d(b,v')=d(b,v)=10\de$. Since geodesic triangles are $\de$-fine, $d(v,v')\le\de$.
For every $x\in supp\,\fff(b,a)$,
\begin{eqnarray*}
d(x,b')&&\le d(x,v)+d(v,v')+d(v',b')\\
&&\le
   8\de+\de+\big[d(b,b')- 10\de\big]\\
&&\le d(b,b')-1,\\
d(a,x) && \le d(a,v)+d(v,x)\\
&&\le
  \big[d(a,b)-10\de\big]+ 8\de\\
&& \le  d(a,b)-1. 
\end{eqnarray*}
Therefore
$$d(a,x)+d(a,b')<d(a,b)+d(a,b').$$ 
Hence the induction hypotheses apply
to the vertices $a$, $x$, and $b'$, giving
\begin{equation}
\label{cond3-2}
\big|r(a,x)-r(a,b')\big| \le d(x,b')+N\le d(b,b')- 1+ N.
\end{equation}
By Proposition~\ref{la'}(1),
$$\fff(b,a)=\sum_{x\in B(v,8\de)} \al_x x$$
for some non-negative coefficients
$\al_x$ summing up to~1. 
By the definition of $r$ and inequality~(\ref{cond3-2}), we have 
\begin{eqnarray*}
&&\big|r(a,b)-r(a,b')\big| \\ 
&&=\Big| r\big(a,\fff(b,a)\big)+ 1- r(a,b')\Big|\\
&&=\left|
     \sum_{x\in B(v,8\de)} \al_x r(a,x)+ 1-
   r(a,b')\right|\\
&&\le
   \left|
     \sum_{x\in B(v,8\de)} \al_x \big[r(a,x)-r(a,b')\big]
   \right|+ 1 \\
&&\le
   \sum_{x\in B(v,8\de)} \al_x \Big|r(a,x)-r(a,b')\Big|+ 1\\
&&\le
   \sum_{x\in B(v,8\de)} \al_x  \big( d(b,b')- 1+ N \big)+ 1\\
&&=  d(b,b')+N.
\end{eqnarray*}

\noi \underline{\bf Case~2.} {\it\ $(a|b')_b\le 10\de$\quad and\quad
$(a|b)_{b'}\le 10\de$.}

\setlength{\unitlength}{1 cm} 
\begin{figure}[h]
  \begin{center}
   \begin{picture}(9,4)
%\graphpaper[1](0,0)(9,4)

\put(0.5,2){\circle*{0.1}}
\put(0,1.95){\small$a$}

\put(8,3.5){\circle*{0.1}}
\put(8.2,3.5){\small$b$}

\put(7,0.5){\circle*{0.1}}
\put(7.2,0.2){\small$b'$}

%from a to b
\qbezier(0.5,2)(6,2)(8,3.5)
\put(1.8,2.4){\footnotesize$p[b,a]$}
%from a to b'
\qbezier(0.5,2)(6,2)(7,0.5)
\put(1.8,1.4){\footnotesize$p[b',a]$}
%from b to b'
\qbezier(8,3.5)(6,2)(7,0.5)

%inscribed triangle
\put(5.93,2.57){\circle*{0.1}}
\put(5.8,2.7){\small$\bar{b'}$}

\put(6.69,1.73){\circle*{0.1}}
\put(5.9,1.31){\circle*{0.1}}
\put(5.7,0.9){\small$\bar{b}$}

\put(7.3,2.5){\tiny $\le$\scriptsize$10\de$}
\put(6.9,1.1){\tiny $\le$\scriptsize$10\de$}

%w
\put(4.4,2.25){\circle*{0.1}}
\put(4.2,2.45){\small$w$}

\put(4.4,1.72){\circle*{0.1}}
\put(4,1.3){\small$v'$}

\put(5.2,2.39){\circle*{0.1}}
\put(5.3,2){\small$v$}

\put(5.6,3.7){\circle*{0.1}}
\put(5.2,3.6){\small$x$}

\put(5,0.5){\circle*{0.1}}
\put(5.2,0.3){\small$x'$}

%from v to x
\qbezier(5.2,2.39)(5.4,3.045)(5.6,3.7)
%from v' to x'
\qbezier(4.4,1.72)(4.7,1.11)(5,0.5)
%from w to w
\qbezier(4.4,1.72)(4.4,1.985)(4.4,2.25)
\put(4.5,1.9){\tiny $\le$\scriptsize$\de$}
   \end{picture}
  \end{center}
 \caption{\label{case2} Proposition~\ref{N}, Case~2.}
\end{figure}
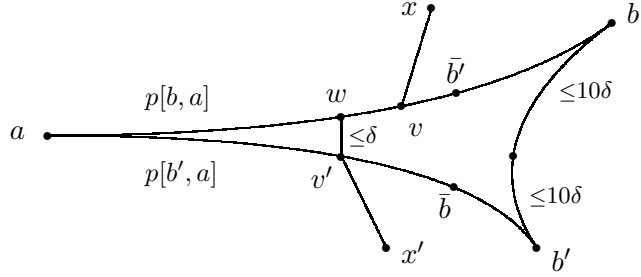

Since $d(a,b)+d(a,b')> 40\de$ and $d(b,b')=(a|b')_b+(a|b)_{b'}\le 20\de$,
we have  $d(a,b)> 10\de$ and $d(a,b')> 10\de$.
Then, by the definition of $r$,
\begin{eqnarray}
\label{cond3-3}
&&\big|
    r(a,b)-r(a,b')
  \big|\\
&&=
  \Big|
    r\big(a,\fff(b,a)\big)+ 1- r\big(a,\fff(b',a)\big)- 1
  \Big| \nonumber \\
&&= \Big|r\big(a,\fff(b,a)- \fff(b',a)\big)\Big|.
\nonumber
\end{eqnarray}
The 0-chain $\fff(b,a)-\fff(b',a)$  can be represented
in the form $f_+-f_-$, where $f_+$ and $f_-$ are 0-chains with
non-negative coefficients and disjoint supports.
By Proposition~\ref{la'}(6),
\begin{eqnarray*}
|f_+|_1+|f_-|_1&&= |f_+ - f_-|_1\\
&&=\big|\fff(b,a)-\fff(b',a)\big|_1 \\
&&\le 2\la'.\end{eqnarray*}
Since the coefficients of the 0-chain $f_+ - f_-=\fff(b,a)-\fff(b',a)$
sum up to 0, then
\begin{equation}
\label{f+-}
|f_+|_1=|f_-|_1\le\la'.
\end{equation}
With the notations $v:=p[b,a](10\de)$, $v':=p[b',a](10\de)$, we have
\begin{eqnarray*}
&&supp\,f_+\se
  supp\,\fff(b,a)\se
  B(v,8\de)\quad\mbox{and}\\
&&supp\,f_-\se
  supp\,\fff(b',a)\se
  B(v',8\de)
\end{eqnarray*}
(see Figure~\ref{case2}).
Since geodesic triangles are $\delta$-fine, there exists a point
$w$ on $p[b,a]$ such that $d(a,w)=d(a,v')$ and $d(w,v')\leq \delta.$
We first assume that $d(a,w)\leq d(a,v).$
We have 
\begin{eqnarray*} d(v,v') && \leq d(v,w)+ d(w,v')\\
&&\leq d(w, \overline{b'})+\delta \\
&& = d(v', \overline{b})+\delta \\
&&\leq   11\delta,
\end{eqnarray*}
where $\overline{b'}$ and $\overline{b}$ are inscribed points as in the definition
of $\delta$-fine triangle in section~\ref{s_hyp}.
If $d(a,w)>d(a,v),$ we can apply the same argument to prove
$d(v,v')\leq 11 \delta$ by interchanging $v'$ with $v$.

Hence by Proposition 2(2), for each $x\in supp\,f_+$ and $x'\in supp\,f_-$,
\begin{eqnarray*}
d(x,x') && \le d(x,v)+d(v,v') +d(v', x')\\ 
&& \le 8\de+ 11\de+  8\de \\
&&= 27\de. 
\end{eqnarray*}

Also
$d(a,x)+d(a,x')< d(a,b)+d(a,b')$, so the induction
hypotheses for the vertices $a$, $x$, and $x'$ apply, giving
\begin{eqnarray}
\label{ineq1}
  \big|r(a,x)-r(a,x')\big|&& \le  d(x,x')+N\\
&&\le  27\de +N
\nonumber
\end{eqnarray}
for each $x\in supp\,f_+$ and $x'\in supp\,f_-$. Then we continue
equality~(\ref{cond3-3}) using (\ref{f+-}), (\ref{ineq1}),
linearity of $r$ in the second variable,
and the definition of $N$ in (\ref{N-def}):
\begin{eqnarray*}
\big|r(a,b)-r(a,b')\big|&& =
  \Big|r\big(a,\fff(b,a)-\fff(b',a)\big)\Big|\\
&&=
  \Big|r\big(a,f_+)-r(a,f_-)\Big|\\
&&\le \la'\cdot \big[27\de+ N\big]\\
&&\le  N\le  d(b,b')+N.
\end{eqnarray*}
Proposition~\ref{N} is proved.\qed

Let $\ep: C_0(G,\qq)\to\qq$ be the augmentation map taking each
0-chain to the sum of its coefficients. A 0-chain~$z$ with
$\ep(z)=0$ is called a {\sf 0-cycle}.

\begin{ppp} 
\label{D} 
There exists a constant $D\ge 0$ such that, for each $a\in G$ and each
0-cycle~$z$,
$$\big| r(a,z) \big|\le D\, |z|_1\, diam \big(supp(z)\big).$$
\end{ppp}
\pf It suffices to consider the case $z=b-b'$, where $b$ and $b'$
are vertices with $d(b,b')=1$. But this case is immediate from
Proposition~\ref{N} by taking $D:=\frac{1}{2}(1+N)$.
\qed

\begin{ttt} 
\label{r}
For a hyperbolic group $G$, the function $r: G\times G\to\qq_{\ge 0}$
from Definition~\ref{def-r} satisfies the following properties.
\begin{itemize}
\item [(1)]$r$ is $G$-equivariant, i.e. $r(a,b)= r(g\cdot a, g\cdot b)$
for $g,a,b\in G$.
\item [(2)] $r$ is Lipschitz equivalent to the distance function. 
      More precisely, we have 
      $$\ds\frac{1}{10\de}\,d(a,b)\le r(a,b)\le d(a,b)$$ for all $a,b\in G$.
\item [(3)] There exist constants $C\ge 0$ and 
$0\le\mu<1$ such that, for all \mbox{$a, a', b, b'\in G$} 
with $d(a,a')\le 1$ and $d(b,b')\le 1$,
$$\big| r(a,b)- r(a',b) - r(a,b')+ r(a',b') \big|\le  C\mu^{d(a,b)}.$$
In particular,
if $d(a,a')\leq 1$ and $d(b,b')\leq 1$, then
$$\big| r(a,b)- r(a',b)- r(a,b')+ r(a',b') \big| \to 0\quad \mbox{as}\quad
d(a,b)\to\infty.$$
\end{itemize}
\end{ttt}
\pf 
\boldmath $(1)$ \unboldmath
The $G$-equivariance of $r$ follows from the definition of $r$
and Proposition~\ref{la'}(4).\\
\boldmath $(2)$ \unboldmath
Using the assumption that $\delta\geq 1$ and the definition of $r$,
the inequalities
\begin{equation*}
\frac{1}{10\de}d(a,b)\le r(a,b)\le d(a,b)
\end{equation*}
can be shown by an easy induction on $d(a,b)$.
The remaining part \boldmath $(3)$ \unboldmath
immediately follows from the following proposition.
\begin{ppp}
\label{AB-rho}
There exist constants $A> 0$, $B> 0$, and $0< \rho <1$ such that,
for all 
$a,a',b,b'\in G$ with $d(a,a')\le 1$ and $d(b,b')\le 30\de$, 
$$\big| r(a,b)- r(a',b)- r(a,b')+ r(a',b') \big|\le 
\big( A\ d(b,b')+ B\big)\,\rho^{d(a,b)+d(a,b')}.$$
\end{ppp}
\pf  Let $D\ge 0$ be the constant from Proposition~\ref{D}, $L\ge 0$ and $0\le\la <1$
be the constants from Propositions~\ref{c-c-c}(5) and~\ref{la'}(5), $\de\ge 1$
be an integral hyperbolicity (fine-triangles) constant, and $\w_7$ be
the number of vertices in a ball of radius $7\de$ in $\G$.

Now we define constants $A$, $B$ and $\rho$. Since the inequality
obviously holds when
$b=b'$, we will assume that $d(b,b')\ge 1$. Then
constant $A>0$ can be chosen large enough so that
\begin{itemize}
\item the desired inequality is satisfied whenever $d(a,b)+d(a,b')\le 100\de$,
     $\rho\ge \sqrt{\la}$, and $B> 0$, and
\item $32D\de L\big(\sqrt{\la}\big)^{-32\de}< A$.
\end{itemize}
 So from now on we can assume that $d(a,b)+d(a,b')>
100\de$. Also the choice of $A$ implies that inequalities
$$1-\frac{A}{Al+B} +
 \frac{32D\de L\big(\sqrt{\la}\big)^{t-32\de}}{(Al+B)\rho^{t-18\de}}\le
1-\frac{A}{Al+B}+
 \frac{32D\de L\big(\sqrt{\la}\big)^{-32\de}}{Al+B}< 1$$
hold for all $B>0$, $\sqrt{\la}\le \rho<1 $, $1\le l\le 30\de$, and $t\ge 0$.
Therefore, we can pick $B>0$ sufficiently large and $\rho <1$ sufficiently close to~1
so that the inequalities
\begin{eqnarray}
\nonumber
&& 1-\frac{A}{Al+B} +
 \frac{32D\de L\big(\sqrt{\la}\big)^{t-32\de}}{(Al+B)\rho^{t-18\de}}\le
\rho^{18\de}\qquad\mbox{and}\\
\nonumber
&& \Big(1- \frac{1}{\w_7} \Big) \frac{30\de A+ B}{B}+
 \frac{64D\de L \big(\sqrt{\la}\big)^{t-32\de}}{B\rho^{t-36\de}}\le \rho^{36\de}
\end{eqnarray}
are satisfied for all $1\le l\le 30\de$ and all $t\ge 0$. The above inequalities
rewrite as
\begin{eqnarray}
\label{rho1}
&& \big( A(l-1)+ B\big)\,\rho^{t-18\de}+
32D\de L\left(\sqrt{\la}\right)^{t-32\de} \le  \big(Al+ B\big)\,\rho^t
\qquad \mbox{and}\\
\label{rho2}
&& \Big(1- \frac{1}{\w_7} \Big) (30\de A+ B)\,\rho^{t-36\de}+
   64D\de L\left(\sqrt{\la}\right)^{t-32\de}\le B\,\rho^t
\end{eqnarray}
and they are satisfied for all $1\le l\le 30\de$ and all $t\ge 0$.

The proof of the proposition proceeds by induction on $d(a,b)+d(a,b')$.
We consider the following two cases.

\noi {\bf \underline{Case~1}.} {\it \ $(a|b)_{b'}> 10\de$\quad or\quad
$(a|b')_b > 10\de$.}

\setlength{\unitlength}{1 cm} 
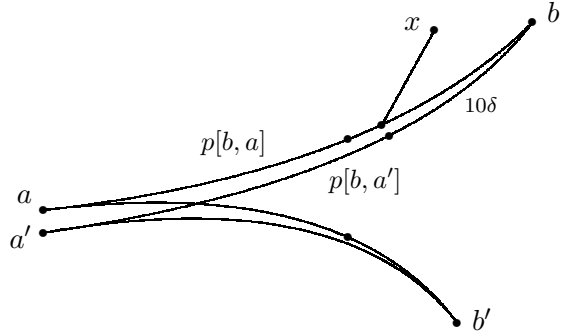
\begin{figure}[h]
  \begin{center}
   \begin{picture}(8,4.5)
%\graphpaper[1](0,0)(8,4.5)

\put(0.5,1.5){\circle*{0.1}}
\put(0.15,1.6){\small$a$}

\put(0.5,1.2){\circle*{0.1}}
\put(0.05,1){\small$a'$}

\put(7,4){\circle*{0.1}}
\put(7.2,4){\small$b$}

\put(6,0){\circle*{0.1}}
\put(6.2,-.1){\small$b'$}

%from a to b
\qbezier(0.5,1.5)(5,2)(7,4)
\put(2.6,2.3){\footnotesize$p[b,a]$}
%from a' to b
\qbezier(0.5,1.2)(5.4,1.9)(7,4)
\put(4.3,1.8){\footnotesize$p[b,a']$}
%from a to b'
\qbezier(0.5,1.5)(4.5,2)(6,0)
%from a' to b'
\qbezier(0.5,1.2)(4.5,1.9)(6,0)

%inscribed triangle
\put(4.55,2.45){\circle*{0.1}}
\put(4.55,1.145){\circle*{0.1}}

\put(6.1,2.8){\tiny$10\de$}

%v
\put(5,2.64){\circle*{0.1}}
%v'
\put(5.1,2.49){\circle*{0.1}}

\put(5.7,3.9){\circle*{0.1}}
\put(5.3,3.9){\small$x$}

%from v to x
\qbezier(5,2.64)(5.35,3.27)(5.7,3.9)
   \end{picture}
  \end{center}
 \caption{\label{case1-AB-rho} Proposition~\ref{AB-rho}, Case~1.}
\end{figure}
Without loss of generality, $(a|b')_b> 10\de$ (interchange $b$
and $b'$ otherwise).
The 0-cycle $f(b,a)- f(b,a')$  can be uniquely represented as $f_+- f_-$,
where $f_+$ and $f_-$ are 0-chains with non-negative coefficients,
disjoint supports, and of the same $\ell^1$-norm. We have
$$f(b,a)=f_0+ f_+ \quad {\rm and} \quad f(b,a')=f_0+ f_-$$ 
for some 0-chain $f_0$ with non-negative coefficients
(actually $f_0= \min \big\{f(b,a),f(b,a')\big\}$).
Denote $\al:= |f_+|_1= |f_-|_1= \ep(f_+)= \ep(f_-)$, where $\ep$ is
the augmentation map.
Since $d(a,a')\le 1$, then 
\begin{eqnarray*}
(a|a')_b 
&&\ge \frac{1}{2}\, \big[ d(a,b)+d(a',b)- 1\big]\\
&&\ge \frac{1}{2}\, \big[ d(a,b)+d(a,b')-32\de\big],  
\end{eqnarray*}
and by Proposition~\ref{c-c-c}(5),
\begin{eqnarray}
\label{alpha}
\al&&= \frac{1}{2}\,\Big| f(b,a)- f(b,a')\Big|_1\\
&&\le
\frac{1}{2}L\la^{(a|a')_b}\nonumber \\ &&\le
\frac{1}{2}L\left(\sqrt{\la}\right)^{d(a,b)+d(a,b')-32\de}.
\nonumber
\end{eqnarray}
By the definition of hyperbolicity in section~\ref{s_hyp}
and the assumption $(a|b')_b>10\de$, we have
$$ d\big(p[b,a](10\de), p[b,a'](10\de)\big) \leq \de.$$ 
Hence there exists  a vertex
$x_0\in B\big(p[b,a](10\de),8\de\big)\cap B\big(p[b,a'](10\de),8\de\big)$.
By the definitions of $r$ and~$\fff$,
\begin{eqnarray*}
&& \big| r(a,b)-  r(a',b)- r(a,b')+ r(a',b') \big|\\
&& = \Big| r\big(a,\fff(b,a)\big)+ 1- r\big(a',\fff(b,a')\big)- 1-
           r(a,b')+ r(a',b') \Big|\\
&& =   \Big| r\big(a,star(f_0+ f_+)\big) - r\big(a',star(f_0+ f_-)\big) -
             r(a,b')+ r(a',b') \Big|\\
&& \le  \Big| r\big(a,star(f_0)+\al x_0\big)-
              r\big(a',star(f_0)+ \al x_0\big) -
              r(a,b')+ r(a',b')  \Big| +\\
&& +  \Big|  r\big(a, star(f_+)- \al x_0\big) \Big| + 
   \Big|  r\big(a',\al x_0- star(f_-)\big) \Big|.\\ 
\end{eqnarray*} 

\noi Now we bound each of the three terms in the last sum.
We number these terms consecutively as $T_1, T_2, T_3$.

\noi  {\it \underline{Term~$T_1$}.} \ Using the same argument 
 as in Case 1 in the proof of Proposition 5, one checks that, for each
$$x\in supp\big( star(f_0)+ \al x_0 \big) \se
B(p[b,a](10\de),8\de)\cap B(p[b,a'](10\de),8\de),$$
the following conditions hold:
\begin{eqnarray*}
&& d(x,b')\le d(b,b')-1\le 30\de  \qquad {\rm and} \\
&& d(a,b)+ d(a,b')- 18\de
\le d(a,x)+d(a,b')
\le d(a,b)+d(a,b')-1.
\end{eqnarray*}
In particular, the induction hypotheses are satisfied for
the vertices $a,a',x,b'$, giving 
\begin{eqnarray*}
&&  \big| r(a,x)- r(a',x)- r(a,b')+  r(a',b')
    \big| \\
&&  \le \big( A\ d(x,b')+ B\big)\,\rho^{d(a,x)+d(a,b')}\\
&&\le
        \Big(A\big(d(b,b')-1\big)+ B\Big)\,\rho^{d(a,b)+d(a,b')-18\de}.
\end{eqnarray*}
Since $star(f_0)+\al x_0$ is a convex combination,
by linearity of $r$ in the second variable,
\begin{eqnarray*}
T_1=&& \Big| r\big(a,star 
(f_0)+\al x_0\big)- r\big(a',star(f_0)+ \al x_0\big)-
             r(a,b')+  r(a',b')
       \Big|\\
\le && \Big(A\, \big(d(b,b')-1\big)+ B\Big)\,\rho^{d(a,b)+d(a,b')-18\de}.
\end{eqnarray*}
 
\noi {\it \underline{Terms~$T_2$ and~$T_3$}.} \ Since
$star(f_+)- \al x_0$ is a 0-cycle supported in a ball of radius~$8\de$,
by Proposition~\ref{D} and inequality~(\ref{alpha}), 
\begin{eqnarray*} 
T_2= && \Big|  r\big(a, star(f_+)- \al x_0\big) \Big|\\
\le && 
           D\ \Big| star(f_+)- \al x_0 \Big|_1\cdot 16\de \\
 \le && D\cdot 2\al\cdot 16\de\\
\le && 
    16D\de  L\left(\sqrt{\la}\right)^{d(a,b)+d(a,b')-32\de}.
\end{eqnarray*} 

Analogously,  
\begin{equation*} 
  T_3=  \Big| r\big(a',\al x_0- star(f_-)\big) \Big|\le
   16D\de  L\left(\sqrt{\la}\right)^{d(a,b)+d(a,b')-32\de}.
\end{equation*}

Combining the three bounds above and using the definition of $B$ and $\rho$
(inequality~(\ref{rho1})),
\begin{eqnarray*}
&& \big| r(a,b)- r(a',b)- r(a,b')+ r(a',b') \big|\\
&&\le T_1+ T_2+ T_3 \\
&& \le \Big( A\,\big(d(b,b')-1\big)+B \Big)\,\rho^{d(a,b)+d(a,b')-18\de}+
       32D\de L\left(\sqrt{\la}\right)^{d(a,b)+d(a,b')-32\de}\\
&& \le \big( A\,d(b,b')+ B\big)\,\rho^{d(a,b)+d(a,b')}.
\end{eqnarray*}
This finishes Case~1.

\noi{\bf \underline{Case~2}.}
{\it \ $(a|b)_{b'}\le 10\de$\quad and\quad $(a|b')_b\le 10\de$.}

\setlength{\unitlength}{1 cm} 
\begin{figure}[h]
  \begin{center}
   \begin{picture}(9,4)
%\graphpaper[1](0,0)(9,4)

\put(0.5,2){\circle*{0.1}}
\put(0,1.95){\small$a$}

\put(0.5,1.7){\circle*{0.1}}
\put(0,1.6){\small$a'$}

\put(8,3.5){\circle*{0.1}}
\put(8.2,3.5){\small$b$}

\put(7,0.5){\circle*{0.1}}
\put(7.2,0.2){\small$b'$}

%from a to b
\qbezier(0.5,2)(6,2)(8,3.5)
%from a to b'
\qbezier(0.5,2)(6,2)(7,0.5)

%from a' to b
\qbezier(0.5,1.7)(6,1.8)(8,3.5)
%from a' to b'
\qbezier(0.5,1.7)(6,1.8)(7,0.5)

%inscribed triangle
\put(5.93,2.57){\circle*{0.1}}
\put(5.9,1.31){\circle*{0.1}}

\put(6.5,3.2){\scriptsize$10\de\!\ge$}
\put(6.6,1.1){\scriptsize$\le\!\! 10\de$}

%w
\put(4.4,2.25){\circle*{0.1}}

%v
\put(5.2,2.39){\circle*{0.1}}
\put(4.95,2.5){\small$v$}

\put(5.25,2.23){\circle*{0.1}}

%v'
\put(4.4,1.72){\circle*{0.1}}
\put(4,1.68){\small$v'$}
\put(4.38,1.53){\circle*{0.1}}

\put(5.6,3.7){\circle*{0.1}}
\put(5.2,3.6){\small$x$}

\put(5,0.5){\circle*{0.1}}
\put(5.2,0.3){\small$x'$}

%from v to x
\qbezier(5.2,2.39)(5.4,3.045)(5.6,3.7)
%from v' to x'
\qbezier(4.4,1.72)(4.7,1.11)(5,0.5)

%from v' to w
\qbezier(4.4,1.72)(4.4,2)(4.4,2.25)
\put(4.45,1.86){\scriptsize $\le\!\! \de$}

   \end{picture}
  \end{center}
 \caption{\label{case2-AB-rho} Proposition~\ref{AB-rho}, Case~2.}
\end{figure}
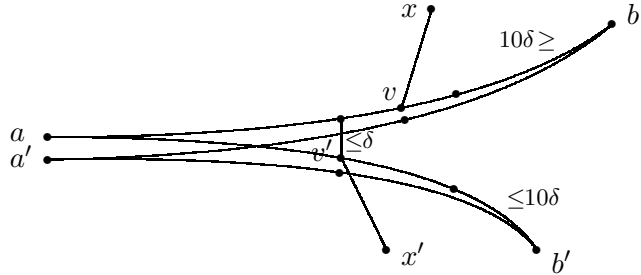

As in Case~1, we have
\begin{eqnarray}
\nonumber & f(b,a)- f(b,a')= f_+- f_-,\qquad
            f(b,a)=f_0+ f_+,\quad  f(b,a')= f_0+ f_-,  &\\
\nonumber & \al:= |f_+|_1= |f_-|_1=\ep(f_+)= \ep(f_-), &\\
\nonumber &
\ds\al\le \frac{1}{2}L\la^{(a|a')_b}\le
\frac{1}{2}L\left(\sqrt{\la}\right)^{d(a,b)+d(a,b')-32\de},&
\end{eqnarray}
where $f_+$, $f_-$, and $f_0$ are 0-chains with non-negative coefficients, 
and $f_+$ and $f_-$ have disjoint supports.
Analogously, interchanging $b$ and $b'$, 
\begin{eqnarray*}
& f(b',a)- f(b',a')= f'_+- f'_-,\qquad
  f(b',a)=f'_0+ f'_+,\quad  f(b',a')=f'_0+ f'_-, &\\
&  \al':= |f'_+|_1= |f'_-|_1=\ep(f'_+)= \ep(f'_-), &\\
&
\ds\al'\le \frac{1}{2}L\la^{(a|a')_{b}}\le
\frac{1}{2}L\left(\sqrt{\la}\right)^{d(a,b)+d(a,b')-32\de},&
\end{eqnarray*}
where  $f'_+$, $f'_-$, and  $f'_0$ are 0-chains with non-negative
coefficients, and $f'_+$ and $f'_-$ have disjoint supports.

Denote $v:=p[b,a](10\de)$ and $v':=p[b',a](10\de)$.
By the conditions of Case~2 and $\de$-hyperbolicity of $\G$,
using the same argument as in Case 2 in the proof of Proposition~\ref{N}, we
obtain
$d(v, v') \leq 11\de$.
Let $x_0$ be a vertex closest to the mid-point of
a geodesic path connecting $v$ to $v'$. 
Proposition~\ref{c-c-c}(2) implies that 
\begin{eqnarray*}
\label{x0}
&& supp\ f_0\cup supp\ f'_0\se B(x_0, 7\de)\qquad\mbox{and}\\
&& supp\ f_-\cup supp\ f_+\cup  supp\ f'_-\cup
supp\ f'_+ \se B(x_0, 8\de).
\end{eqnarray*}
By the definition of $r$,
\begin{eqnarray*}
&& \big| r(a,b)- r(a',b)- 
         r(a,b')+ r(a',b') \big|\\
&&  =  \Big| r\big(a,\fff(b,a)\big)- r\big(a',\fff(b,a')\big)- 
             r\big(a,\fff(b',a)\big)+ r\big(a',\fff(b',a')\big) \Big|\\
&& \le   \Big| r\big(a,star(f_0+ f_+)\big)-
               r\big(a',star(f_0+ f_-)\big)-
               r\big(a,star(f'_0+ f'_+)\big)+
               r\big(a',star(f'_0+ f'_-)\big)
       \Big|\\
&& \le  \Big| r\big(a, star(f_0)+  \al  x_0 -
                        star(f'_0)- \al' x_0 \big)-
              r\big(a', star(f_0)+  \al  x_0 -
                         star(f'_0)- \al' x_0 \big)\Big| +\\
&& +  \Big|  r\big(a, star(f_+)\big)- r\big(a,\al x_0\big)   \Big| + 
      \Big|  r\big(a',\al x_0\big)- r\big(a',star(f_-)\big)  \Big|+\\
&&+   \Big|  r\big(a,\al' x_0\big)- r\big(a, star(f'_+)\big) \Big| + 
      \Big|  r\big(a',star(f_-')\big)- r\big(a',\al' x_0\big)  \Big|.\\
\end{eqnarray*}
\noi Now we bound each of the five terms in the last sum.
We number these terms consecutively as $S_1$, ...,  $S_5$.

\noi  {\it \underline{Term~$S_1$}.}
One checks that, for each
\begin{eqnarray*}
&& x\in supp\big( star(f_0)+ \al x_0 \big) \se
B(v,8\de)\cap B(p[b,a'](10\de),8\de)\quad \mbox{and}\\
&& x'\in supp\big( star(f'_0)+ \al' x_0 \big) \se
B(v',8\de)\cap B(p[b',a'](10\de),8\de),
\end{eqnarray*}
the following conditions hold:
\begin{eqnarray*}
&&  d(x,x')\le 30\de \quad \mbox{and}\\
&&  d(a,b)+ d(a,b')- 36\de  \le d(a,x)+d(a,x')
\le d(a,b)+d(a,b')-1.
\end{eqnarray*}
In particular, the induction hypotheses are satisfied for
the vertices $a,a',x,x'$, giving 
\begin{eqnarray}
\label{ind2}
&&  \big| r(a,x)- r(a',x)- r(a,x')+ r(a',x') \big| \\
\nonumber &&  \le \Big( A\ d(x,x')+ B\Big)\,\rho^{d(a,x)+d(a,x')}\\ 
\nonumber && \le
        \Big(30\de A+ B\Big)\,\rho^{d(a,b)+d(a,b')-36\de}.
\end{eqnarray}

Recall that $\w_7$ is the number of vertices in a ball of radius $7\de$.
Let $\be$ be the (positive) coefficient of $x_0$ in the 0-chain $star(f_0)$,
and $\be'$ be the (positive) coefficient of $x_0$ in the 0-chain $star(f'_0)$.
Without loss of generality, we can assume
 $\big|star(f_0)\big|_1\le \big|star(f'_0)\big|_1$.
Since $x_0$ was chosen so that $supp\ f_0\cup supp\ f'_0\se B(x_0,7\de)$,
by the definition of $star$,
we have 
\begin{eqnarray*}
&&\be= {\frac{1}{\w_7}}\, \big|f_0\big|_1=
  {\frac{1}{\w_7}}\, \big|star(f_0)\big|_1 \le
  {\frac{1}{\w_7}}\, \big|star(f'_0)\big|_1=
  {\frac{1}{\w_7}}\, \big|f'_0\big|_1= \be'\qquad\mbox{and}\\
&&\al -\al'=
\Big( 1- \big|f_0 \big|_1\Big)-
\Big( 1- \big|f'_0 \big|_1\Big)=
\big|f'_0 \big|_1- \big|f_0 \big|_1= \w_7(\be'-\be)\ge 0.
\end{eqnarray*}
Therefore,
\begin{eqnarray*}
&& \Big| star(f_0)+  \al  x_0-
         star(f'_0)- \al' x_0  \Big|_1\\
&& \le \Big|  star(f_0)- \be x_0  \Big|_1 +
       \Big|  \be' x_0- star(f'_0)  \Big|_1 +
       \Big|  \big[ (\al- \al')- (\be' -\be) \big] x_0  \Big|_1\\
&& =   \Big( \big|  star(f_0) \big|_1 - \be \Big)+
       \Big( \big|  star(f'_0) \big|_1 - \be' \Big)+
       \big( \be' -\be \big) \big( \w_7-1 \big)\\
&& =   \Big( \big| f_0  \big|_1 - \be \Big)+
       \Big( \big| f'_0 \big|_1 - \be' \Big)+
       \Big( \big|f'_0 \big|_1- \big|f_0 \big|_1 \Big)
        \Big(1- \frac{1}{\w_7} \Big)\\
&& =   \big| f_0  \big|_1  \Big(1- \frac{1}{\w_7} \Big) +
       \big| f'_0 \big|_1  \Big(1- \frac{1}{\w_7} \Big) +
       \Big( \big|f'_0 \big|_1- \big|f_0 \big|_1 \Big)
        \Big(1- \frac{1}{\w_7} \Big)\\
&& =   2 \big| f'_0 \big|_1   \Big(1- \frac{1}{\w_7} \Big)\\ &&
 \le
       2 \Big(1- \frac{1}{\w_7} \Big).
\end{eqnarray*}

Since $\big[ star(f_0)+  \al  x_0 \big]- \big[ star(f'_0)+ \al' x_0 \big]$
 is a 0-cycle, it is of the
 form $h_+- h_-$, where $h_+$ and $h_-$ are 0-chains with
non-negative coefficients, disjoint supports and of the same $\ell^1$-norm,
so we can define 
$$\g := |h_+|_1= |h_-|_1= \ep(h_+)= \ep(h_-).$$
By the above inequality,
$$\g= {\frac{1}{2}}\ |h_+- h_-|_1\le  1- \frac{1}{\w_7},$$
then, by~(\ref{ind2}) and linearity of $r$
in the second variable,
\begin{eqnarray*}
S_1\ = && \Big| r(a, h_+- h_-)-
                r(a',h_+- h_-)
          \Big|\\
= && \Big| r(a,h_+)- r(a',h_+)-
           r(a,h_-)+ r(a',h_-)
     \Big|\\
\le &&  \g\cdot \Big(30\de A+ B\Big)\,\rho^{d(a,b)+d(a,b')-36\de}\\
\le
&&        \Big(1- \frac{1}{\w_7} \Big)
        \Big(30\de A+ B\Big)\,\rho^{d(a,b)+d(a,b')-36\de}.
\end{eqnarray*}
\noi  {\it \underline{Terms~$S_2-S_5$}.} \ Analogously to term~$T_2$
in Case~1,
\begin{eqnarray*} 
& S_2\ \le  16D\de  L\left(\sqrt{\la}\right)^{d(a,b)+d(a,b')-32\de},\qquad
  S_3\ \le  16D\de  L\left(\sqrt{\la}\right)^{d(a,b)+d(a,b')-32\de},&\\
& S_4\ \le  16D\de  L\left(\sqrt{\la}\right)^{d(a,b)+d(a,b')-32\de},\qquad
  S_5\ \le  16D\de  L\left(\sqrt{\la}\right)^{d(a,b)+d(a,b')-32\de}.&\\
\end{eqnarray*}
Combining the bounds for the five terms above and using the definition
of $B$ and $\rho$ (inequality~(\ref{rho2})),
\begin{eqnarray*}
&& \big| r(a,b)- r(a',b)-
          r(a,b')+ r(a',b') \big|\\
&& \le S_1+ S_2+ S_3+ S_4+ S_5\\
&& \le \Big(1- \frac{1}{\w_7} \Big)
       (30\de A+ B)\,\rho^{d(a,b)+d(a,b')-36\de}+
   64D\de L\left(\sqrt{\la}\right)^{d(a,b)+d(a,b')-32\de}\\
&& \le B\,\rho^{d(a,b)+d(a,b')}\\
&&\le
   \big( A\,d(b,b')+ B\big)\,
\rho^{d(a,b)+d(a,b')}.
\end{eqnarray*} 
Proposition~\ref{AB-rho} and Theorem~\ref{r} are proved.
\qed

\section{More properties of $r$.}

In this section, we prove two distance-like inequalities for
the function $r$ introduced in the previous section.

As before, let $G$ be a hyperbolic group and 
$\Gamma$ be the Cayley graph of $G$ with respect to a
finite generating set.  For any subset $A\subseteq \Gamma$, denote
$$N_G(A,R) :=\{ x\in G\ \big|\  d(x, A)\leq R\}.$$

\begin{ppp}
\label{C1}
There exists $C_1\geq 0$ with the following property.
If $a,b\in G$, $\g$ is a geodesic in~$\G$ connecting  $a$ and $b$,
$x\in G\cap\g$, $\g'$ is the part of $\g$ between $x$ and $b$, and
$c\in N_G (\g', 9\delta)$, then
$$\big|r(a,c)-r(a,x)-r(x,c)\big| \leq C_1 \qquad (Figure~\ref{fig-axc}).$$
\end{ppp}
\setlength{\unitlength}{.7cm} 
\begin{figure}[h]
  \begin{center}
   \begin{picture}(12,4.2)
%\graphpaper[1](0,0)(12,4.2)
\put(1,1.5){\circle*{0.15}}
\put(.3,1.4){\small$a$}

\put(11,1.5){\circle*{0.15}}
\put(11.4,1.35){\small$b$}

\put(8,2.5){\circle*{0.15}}
\put(8.4,2.35){\small$c$}

\put(8,1.5){\circle*{0.15}}
\put(8.3,1){\small$c'$}

\put(4,1.5){\circle*{0.15}}
\put(4.3,1){\small$x$}

\put(6,0.5){\footnotesize$\g$}

%from a to b
\qbezier(1,1.5)(6,1.5)(11,1.5)

%from x to c
\qbezier(4,1.5)(6,2)(8,2.5)

%from a to c
\qbezier(1,1.5)(4.5,1.7)(8,2.5)
   \end{picture}
  \end{center}
 \caption{\label{fig-axc} Proposition~\ref{C1}.}
\end{figure}
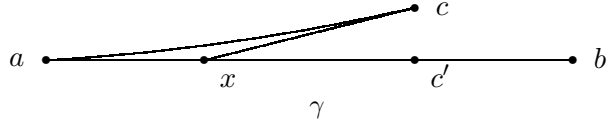

\pf
 Let 
$$C_1 :=(80\de+N+36\de DL)\sum_{k=0}^{\infty} \lambda^{k-18\delta},$$
where $L\geq 1$ and $0< \lambda <1$ are as in Propositions~\ref{c-c-c}(5)
and~\ref{la'}(5), $N$ is as in Proposition~\ref{N}, and $D$ is as
in Proposition~\ref{D}. 
It suffices to show the inequality
$$\big|r(a,c)- r(a,x)-r(x,c)\big| \leq 
(80\de+N+36\de DL)\sum_{k=0}^{d(x,c)} \lambda^{k-18\delta}.$$
We will prove it by induction on $d(x,c)$.

If $d(x,c)\leq 40\delta$, by Proposition 4 and Theorem 6(2) we have
 \begin{eqnarray*} 
&&\big|r(a,c)-r(a,x)-r(x,c)\big| 
    \leq \big|r(a,c)-r(a,x)\big| +r(x,c)\\ 
&& \leq \big(d(c,x)+N \big)+ d(x,c)   \leq 80\de+ N\\
&& \leq (80\de+ N+36\de DL) \sum_{k=0}^{d(x,c)} \la^{k-18\de}.
\end{eqnarray*}
Now we assume that $d(x,c)> 40 \delta.$
There exists a vertex $c'\in \g'$ with $d(c',c)\leq 9\delta$, so
$$d(a,c)\geq d(a,c')-9\delta 
   \geq d(x,c')-9\delta 
   \geq d(x,c)-18\delta 
   >10\delta.$$
Hence by the definition of the function $r$, we have
$$r(a,c)=r\big(a,\bar{f}(c,a)\big)+1\qquad\mbox{and}\qquad
r(x,c)=r\big(x,\bar{f}(c,x)\big)+1.$$
Also
\begin{eqnarray*}
(a|x)_c 
&&= \frac{1}{2}\big[d(c,a)+d(c,x)-d(a,x)\big]\\ 
&&\ge\frac{1}{2} \big[d(c',a)- 9\de+ d(c',x)- 9\delta- d(a,x)\big]\\ 
&&=d(x,c')-9\delta \\
&&\geq d(x,c)-18\delta.
\end{eqnarray*}
By Proposition~\ref{la'}(5),
$$\Big|\bar{f}(c,x)-\bar{f}(c,a)\Big|_1\leq
L \lambda^{(a|x)_c}\leq 
L\lambda^{d(x,c)-18\delta}.$$
This, together with Proposition~\ref{D} and Proposition 2(2), implies that
\begin{eqnarray*}
\Big| r\big(a,\bar{f}(c,a)\big)- r\big(a,\bar{f}(c,x)\big)\Big|
&&= \Big| r\big(a,\bar{f}(c,a)- \bar{f}(c,x)\big)\Big|\\
&&\leq 
DL \lambda^{d(x,c)-18\delta} \, diam\big(supp (\fff(c,a)-\fff(c,x))\big)\\
&& \le 36\de DL \la^{d(x,c)-18\de}.
\end{eqnarray*}

By Proposition~\ref{la'}(2) and~\ref{la'}(7),
for every $y\in supp(\bar{f}(c,x))$, we have 
$$ d(x,y)\leq d(x,c)-1\qquad\mbox{and}\qquad y\in N_G (\g', 9\delta).$$

Hence by the induction hypotheses, we obtain 
\begin{eqnarray*}     
&&\big| r(a,c)-r(a,x)-r(x,c)\big|\\
 &&=
 \Big| \big(r\big(a,\bar{f}(c,a)\big)+1\big)-r(a,x)-\big(r(x,\bar{f}(c,x))+1\big)\Big|\\ 
&&\leq \Big|r\big(a,\bar{f}(c,a)\big)-r\big(a, \bar{f}(c,x)\big)\Big| 
+\Big|r\big(a,\bar{f}(c,x)\big)- r(a,x)- r\big(x,\bar{f}(c,x)\big)\Big|\\ 
&&\leq 36\de DL\la^{d(x,c)-18\de}+
 (80\de+ N+ 36\de DL) \sum_{k=0}^{d(x,c)-1}\lambda^{k-18\delta}\\
&&\le (80\de+ N+ 36\de DL) \sum_{k=0}^{d(x,c)}\lambda^{k-18\delta}.
\end{eqnarray*}
\qed

\begin{ppp} 
\label{M'}
There exists $M'\geq 0$ such that
$$\big|r(a,b)-r(a',b)\big| \leq M'\,d(a,a')$$
for all $a, a', b \in G$.
\end{ppp} 
 
\pf
Recall that $\delta\geq 1$. Let
$$M':=(20\delta+ 3+ 36\de DL) \sum_{k=0}^{\infty} \lambda^{k-19\de}.$$
The Cayley graph $\Gamma$ is a geodesic metric space, hence
it suffices to show the inequality
$\big|r(a,b)-r(a',b)\big|\le M'$ in the case when $d(a,a')=1$.
For that, it suffices to show the inequality
$$\big|r(a,b)-r(a',b)\big| \leq (20\delta+3+36\de
DL)\sum_{k=0}^{d(b,a)}\lambda^{k-19\de}
$$ when $d(a,a')=1$. We will prove it by induction on $d(a,b)$.

If $d(a,b)\leq 10\delta +1$, then by Theorem 6(2) we have
\begin{eqnarray*}
\big|r(a,b)-r(a',b)\big| && \leq r(a,b)+ r(a',b)\\
&& \leq d(a,b)+ d(a',b)\\
&& \leq 20\delta +3 \\
&&\le
(20\delta+ 3+ 36\de DL)\sum_{k=0}^{d(b,a)} \la^{k-19\de}.
\end{eqnarray*}

If $d(a,b)>10\delta +1$, then $d(a',b)>10\delta.$

For every $y\in supp(\fff(b,a))\cup supp(\fff(b,a'))$,
by Proposition 2(2) we have 
$$(a|a')_y=\frac{1}{2}\big[d(y,a)+d(y,a')-d(a,a')\big]\ge
d(b,a)-19\de.$$
Hence by the definition of the function $r$, the induction hypothesis
and Propositions~\ref{la'}(5) and ~\ref{D}, we obtain
\begin{eqnarray*}
&& \big|r(a,b)-r(a',b)\big|\\
 &&= \Big|\big(r(a, \bar{f}(b,a))+1\big)- \big(r(a', \bar{f}(b,a'))+1\big)\Big|\\
&&\leq \big|r(a,\bar{f}(b, a))-r(a',
\bar{f}(b,a))\big|+\big|r(a',\bar{f}(b,a))-r(a',\bar{f}(b,a'))\big|\\
&&\leq  (20\delta+ 3+ 36\de DL)\sum_{k=0}^{d(b,a)-1} \la^{k-19\de}+
  DL\la^{d(b,a)-19\de} diam\Big(supp\big (\fff(b,a)-\fff(b,a')\big)\Big)\\
&&\leq  (20\delta+ 3+ 36\de DL)\sum_{k=0}^{d(b,a)-1} \la^{k-19\de}+
   36\de DL\la^{d(b,a)-19\de}\\
&&\le (20\delta+ 3+ 36\de DL)\sum_{k=0}^{d(b,a)} \la^{k-19\de}.
\end{eqnarray*}   
\qed

\section{Definition and properties of a new metric $\hat{d}$.}

In this section, we use the function $r$ defined in section 3
to construct a $G$-invariant metric $\hat{d}$ on a hyperbolic group $G$
such that $\hat{d}$ is quasi-isometric to the word metric  and prove that
$(G, \hat{d})$ is weakly geodesic and
strongly bolic.

We  define 
$$s(a,b):=\frac{1}{2}\, \big[r(a,b)+ r(b,a)\big]$$
for all $ a,b \in G$.

\begin{ppp}
\label{uvw}
The above function $s$ satisfies the following conditions.
\begin{itemize} 
\item [(a)] There exists $M\ge 0$ such that
$$\big|s(u,v)-s(u,v')\big|\le M\, d(v,v')\quad\mbox{and}\quad
   \big|s(u,v)-s(u',v)\big|\le M\, d(u,u')$$
for all $u, u', v,v'\in G$.
\item [(b)] There exists $C_1\ge 0$ such that if a vertex $w$ lies on a geodesic
 connecting vertices~$u$ and~$v$, then
  $$\big| s(u,v)-s(u,w)-s(w,v)\big|\le C_1.$$
\end{itemize}
\end{ppp}
\pf
{\bf (a)}
Since $s$ is symmetric, it suffices to show only the first inequality.
Since the Cayley graph  $\Gamma$ is a geodesic metric space,
 it suffices to consider only the case  $d(v,v')=1$.
This case follows from Propositions~\ref{N} and~\ref{M'}.

{\bf (b)}
follows from Proposition~\ref{C1}.
\qed

\begin{ppp}
\label{C2}
There exists $C_2\geq 0$ such that
$$s(a,b)\leq s(a,c)+ s(c,b)+C_2$$
for all $ a, b, c\in G$.
\end{ppp}

\pf
Let $\bar{a} \in p[b,c], \bar{c} \in p[a,b], \bar{b}\in p[a,c]$ such 
that 
$$d(b, \bar{c})= d(b, \bar{a}), \quad d(c, \bar{a})= d(c, \bar{b}),
 \quad d(a, \bar{c})= d(a, \bar{b}).$$

By the definition of hyperbolicity, we have
$$ d(\bar{a}, \bar{b})\leq \delta, \quad d(\bar{a}, \bar{c})\leq \delta,
\quad d(\bar{b}, \bar{c})\leq \delta.$$
By Proposition~\ref{uvw},
\begin{eqnarray*} 
 s(a,b)
&& \leq s(a,\bar{c}) + s(\bar{c}, b) +C_1\\ 
&&\leq \big(s(a,\bar{b})+ M\,d(\bar{b},\bar{c})\big)+
    \big( s(\bar{a},b) + M\,d(\bar{c},\bar{a}) \big)+ C_1\\
&&\le s(a,\bar{b})+s(\bar{a},b)+ 2\de M + C_1\\
&&\le \big(s(a,\bar{b})+ s(\bar{b},c)\big)
   +\big(s(c,\bar{a})+ s(\bar{a},b)\big)+ 2\de M+ C_1\\
&&\le s(a,c)+ s(c,b)+ 2\de M+ 3C_1,
\end{eqnarray*} 
so we set $C_2:= 2\de M+ 3C_1$.
\qed

For every  pair of elements $a, b\in G$, we define

\[ \hat{d}(a,b):=\left\{ \begin{array}{ll}
s(a, b) +C_2 & \mbox{if $a\neq b$,} \\
0 & \mbox{if $a=b$.}
\end{array}
\right. \]

\begin{ppp} The function $\hat{d}$ defined above is a metric on $G$.
\end{ppp} 

\pf
By definition, $\hat{d}$ is symmetric, and $\hat{d}(a,b)=0$ iff $a=b$.
The triangle inequality is a direct consequence of Proposition~\ref{C2}.
\qed

\begin{ppp}
\label{d0}
There exist constants $C\ge 0$ and 
$0\le\mu<1$ with the following property.
For all $R\ge 0$ and all \mbox{$a, a', b, b'\in G$}
with $d(a,a')\le R$ and $d(b,b')\le R$,
$$\big| \hat{d}(a,b)- \hat{d}(a',b) - \hat{d}(a,b')+ \hat{d}(a',b') \big|\le 
R^2 C\mu^{d(a,b)-2R}.$$
In particular, if  $d(a,a')\le R$ and $d(b,b')\le R$, then
$$\hat{d}(a,b)-\hat{d}(a', b)- \hat{d}(a,b')+\hat{d}(a',b') \rightarrow 0
   \quad\mbox{as}\quad d(a,b)\rightarrow \infty.$$
\end{ppp}

\pf
Take $C$ and $\mu$ as in Theorem~\ref{r}(3).
Increasing $C$ if needed we can assume that
$a\neq b, a\neq b', a'\neq b, a'\neq b'$.

If $ a=a'$ or $b=b'$, then
$$\hat{d}(a,b)-\hat{d}(a', b)- \hat{d}(a,b')+\hat{d}(a',b')=0.$$

If $d(a,a')=1$ and $d(b,b')=1$, then by Theorem~\ref{r}(3),
\begin{eqnarray*}
&&\big|\hat{d}(a,b)-\hat{d}(a', b)- \hat{d}(a,b')+\hat{d}(a',b')\big|\\
&&= \big|s(a,b)-s(a',b)-s(a,b')+s(a',b')\big|\\
&&\le C\mu^{d(a,b)}. 
\end{eqnarray*}

Without loss of generality, we can assume that  $R$ is  an integer.
In the general  case
$$d(a,a')\le R\qquad \mbox{and}\qquad d(b,b')\le R,$$
pick vertices $a=a_0, a_1, ... , a_R=a'$
with $d(a_{i-1}, a_i)\le 1$ and $b=b_0, b_1, ... , b_R=b'$  with $d(b_{j-1},b_j)\le 1$ and
note that
$d(a_i,b_j)\ge d(a,b)-2R$. Then we have
\begin{eqnarray*}
&&\big|\hat{d}(a,b)-\hat{d}(a', b)- \hat{d}(a,b')+\hat{d}(a',b')\big|\\
&&= \big|s(a,b)-s(a',b)-s(a,b')+s(a',b')\big| \\
&&= \Big|\sum_{i=1}^{R} \sum_{j=1}^{R}
     \big( s(a_{i-1},b_{j-1})-s(a_i,b_{j-1})-s(a_{i-1},b_j)+s(a_i,b_j)\big)\Big| \\
&&\le \sum_{i=1}^{R} \sum_{j=1}^{R}
      \big|s(a_{i-1},b_{j-1})-s(a_i,b_{j-1})-s(a_{i-1},b_j)+s(a_i,b_j)\big| \\
&&\le  R^2 C\mu^{d(a,b)-2R}.
\end{eqnarray*}
\qed

Recall that a metric space $(X, d)$ is said to be
{\sf weakly geodesic}~\cite{KS, KS1} 
if there exists $\de_1\geq 0$ such that, for every pair of points $x$ and $y$ in $X$
and  every $t\in [0, d(x,y)]$, there exists a point $a\in X$ 
such that $d(a,x)\leq t+\de_1$ and $d(a,y)\leq d(x,y)-t+\de_1$.

\begin{ppp} 
\label{weakly-geodesic}
The metric space $(G, \hat{d})$ is weakly geodesic.
\end{ppp}

\pf
Let $x, y\in G$ and $z\in G\cap p[x,y]$.
By the definition of $\hat{d}$ and  Proposition~\ref{uvw}(b), we have
$$ \hat{d}(x,z) +\hat{d}(z,y)-\hat{d}(x,y)\leq C_1+2C_2.$$
It follows that 
$$ \hat{d}(x,z)\leq \hat{d}(x,y)+ C_1 + 2C_2,$$
hence the image of the map
$$\hat{d}(x, \cdot): G\cap p[x,y] \rightarrow [0, \infty),$$
is contained in $[0,\hat{d}(x,y) +C_1 +2C_2].$
Also, the image contains 0 and $\hat{d}(x,y)$.

By Proposition~\ref{uvw}(a), we have
$$\big|\hat{d}(x, z')-\hat{d}(x, z)\big|\leq M$$
when $d(z', z)=1$.
This, together with the fact that $p[x,y]$ is a geodesic path,
implies that the image of the map
$$\hat{d}(x, \cdot): G\cap p[x,y] \rightarrow [0, \hat{d}(x,y)+C_1+2C_2]$$
is $M$-dense in $[0, \hat{d}(x,y)]$, i.e.
for every $t\in [0, \hat{d}(x,y)]$, there exists 
$a \in G\cap p[x,y]$ such that
$$\big|\hat{d}(x,a)-t\big| \leq M.$$
It follows that $\hat{d}(x,a)\leq t+ M$, and
by Proposition~\ref{uvw}(b) we also have 
$$\big| \hat{d}(x,y)- \hat{d}(x,a)-\hat{d}(a,y)\big| \leq C_1 +2C_2.$$
This implies that
\begin{eqnarray*}
\hat{d}(a,y) && \leq \hat{d}(x,y)-\hat{d}(x,a) +C_1 +2C_2\\ 
&&\leq \hat{d}(x,y) -t +M+C_1+2C_2.
\end{eqnarray*}
Therefore $(G,\hat{d})$ is weakly geodesic for
$\de_1:=M+C_1+2C_2$.
\qed

Kasparov and Skandalis introduced the concept of bolicity
 in~\cite{KS, KS1}. 
\begin{ddd}
A metric space $(X,d)$ is said to be {\sf bolic}  if there exists
$\de_2\ge 0$ with the following properties:
\begin{itemize}
\item [(B1)] for any $R>0$, there exists $R'>0$ such that for all  $a, a',b,b'
\in X$ satisfying
$$d(a,a')+d(b,b')\leq R\qquad\mbox{and}\qquad d(a,b)+d(a',b')\geq R',$$
we have $$d(a,b')+d(a',b)\leq d(a, b)+ d(a', b') +2\de_2;\qquad\mbox{and}$$  
\item [(B2)] there exists a map $m: X\times X \rightarrow X$,
such that, for all $x,y, z\in X$, we have
$$2d(m(x,y), z)\leq \big(2d(x,z)^2 +2d(y,z)^2-d(x,y)^2\big)
^{\frac{1}{2}} +4\de_2.$$   
\end{itemize}
$(X,d)$ is called {\sf strongly bolic} if it is bolic and
the above condition (B1) holds for every $\de_2> 0$~\cite{L3}.
\end{ddd}

\begin{ppp}
\label{bolic}
The metric space $(G, \hat{d})$ is strongly bolic.
\end{ppp}

\pf
Proposition~\ref{d0} yields condition~(B1) for all $\de_2>0$.
It remains to show that there exist $\de_2\ge 0$ and
a map $m: G\times G \rightarrow G$, such that, for all $x,y, z\in G$, we have
$$2\hat{d}(m(x,y), z)\leq
\big(2\hat{d}(x,z)^2 +2 \hat{d}(y,z)^2-\hat{d}(x,y)^2\big)^{\frac{1}{2}} +4\de_2.$$

By Proposition~\ref{weakly-geodesic} and its proof, there exists a vertex
$m(x,y)\in G\cap p[x,y]$ such that
\begin{equation}
\label{mxy}
\Big|\hat{d} (x, m(x,y))-\frac{\hat{d}(x,y)}{2}\Big|\leq \delta_1\quad\mbox{and}
 \quad\Big|\hat{d}(m(x,y),y)-\frac{\hat{d}(x,y)}{2} \Big| \leq \delta_1.
\end{equation}
By the definition of $\delta$-hyperbolicity, we know that
either
\begin{itemize}
\item [(1)] there exists $a\in G\cap p[z,y]$ such that 
$d(m(x,y), a)\leq \delta +1$, or
\item [(2)] there exists $b\in G \cap p[x,z]$ such that 
$d(m(x,y), b)\leq \delta +1.$
\end{itemize}

In case (1), we have 
$$\big|\hat{d}(z, m(x,y))-\hat{d}(z,a)\big| \leq \hat{d}(m(x,y), a)\leq \delta +1+C_1,$$
$$\big|\hat{d}(y, m(x,y))-\hat{d}(y,a)\big| \leq \hat{d}(m(x,y),a)\leq \delta +1+C_1.
$$
Hence, by Proposition~\ref{uvw}(b), we obtain 
\begin{eqnarray*}
\hat{d}(z, m(x,y))+\hat{d}(x,y) &&\leq
 \hat{d}(z,a)+\delta +1+C_1 +\hat{d}(x,y)\\ 
&&\le \hat{d}(z,a)+ \delta +1+C_1 +\hat{d}(x, m(x,y))+ \hat{d}(m(x,y),y)\\ 
&&\leq   \hat{d}(z,a)+\hat{d}(a,y)+\hat{d}(x, m(x,y))+ 2\delta +2C_1+2\\
&&\leq \hat{d}(y,z)+\hat{d}(x, m(x,y))+\delta',
\end{eqnarray*}
where $\delta':=2\delta +3C_1 +2$. In case (2), we similarly have
$$ \hat{d}(z, m(x,y))+ \hat{d}(x,y) \leq 
\hat{d}(x,z)+\hat{d}(m(x,y),y) +\delta'.$$
It follows from~(\ref{mxy}) that 
\begin{eqnarray*}
 \hat{d}(z, m(x,y))+\hat{d}(x,y) && \leq 
\sup\big\{\hat{d}(x,z)+\hat{d}(y, m(x,y)),\,  \hat{d}(y,z) +\hat{d}(x, m(x,y))\big\}
+\delta'\\ 
&&\leq \sup \big\{\hat{d} (x,z), \hat{d}(y,z)\big\} +\frac{\hat{d}(x,y)}{2}
 +\delta_1 +\delta'.
\end{eqnarray*}
Hence 
\begin{equation}
\label{2d}
2\hat{d}(z, m(x,y)) \leq 2 \sup \big\{\hat{d}(x,z), \hat{d}(y,z)\big\}-
\hat{d}(x,y)+ 4\delta_2,
\end{equation}
where $\delta_2 :=\frac{\delta_1 +\delta'}{2}.$

If $t, u,$ and $v$ are non-negative real numbers such that 
$|u-v|\leq t$, then
$$(2u-v)^2 \leq 2u^2+2t^2 -v^2.$$
Setting
$t:=\inf \{\hat{d}(x,z), \hat{d}(y,z)\},
u:=\sup\{\hat{d}(x,z), \hat{d}(y,z)\},
v:=\hat{d}(x,y),$
we obtain
$$2 \sup\big\{\hat{d}(x,z), \hat{d}(y,z)\big\}-\hat{d}(x,y)
\leq \big(2\hat{d}(x,z)^2+2 \hat{d}(y,z)^2 -\hat{d}(x,y)^2\big)^{\frac{1}{2}}.$$
Therefore, by~(\ref{2d}),
$$2\hat{d}(z, m(x,y))\leq
\big(2\hat{d}(x,z)^2 +2 \hat{d}(y,z)^2-\hat{d}(x,y)^2\big)^{\frac{1}{2}} +4\delta_2.$$
\qed 

We summarize the results of this section.
\begin{ttt}
\label{metric}
Every hyperbolic group $G$ admits a metric $\hat{d}$ with the following properties.
\begin{itemize}
\item [(1)] $\hat{d}$ is $G$-invariant, i.e.
$\hat{d}(g\cdot x,g\cdot y)=\hat{d}(x,y)$ for all $x,y,g\in G$.
\item [(2)] $\hat{d}$ is quasiisometric to the word metric $d$, i.e.
there exist
$A>0$ and $B\geq 0$ such that
$$\frac{1}{A} \hat{d}(x,y)-B \leq d(x,y) \leq A \hat{d}(x,y) +B$$  
for all $x,y\in G$.
\item [(3)] The metric space $(G,\hat{d})$ is weakly geodesic and 
strongly bolic.
\end{itemize}
\end{ttt}

\section{The Baum-Connes conjecture for hyperbolic groups.}

In this section, we combine Theorem 17 with Lafforgue's work to prove
the main result of this paper.
\begin{ddd}
An action of a topological group $G$ on a topological space $X$ is called {\sf proper}
if the map $G\times X\to X\times X$ given by $(g,x)\mapsto (x,gx)$
is a proper map, that is the preimages of compact subsets are compact.
\end{ddd}

When $G$ is discrete, an action is proper iff it is properly discontinuous,
i.e. if the set $\{g\in G\ \big|\ K\cap gK \not=\emptyset\}$ is finite
for any compact $K\se X$.

The following deep theorem was proved by Lafforgue using Banach
KK-theory.

\begin{ttt}[Lafforgue~\cite{L3}] 
If a discrete group $G$ has property RD, and G  acts properly and
isometrically on a strongly bolic, weakly geodesic,
 and uniformly locally finite metric
space, then the Baum-Connes conjecture holds for $G$.
\end{ttt}

\begin{ttt} 
\label{bc}
The Baum-Connes conjecture holds for 
  hyperbolic groups and their subgroups.
\end{ttt}

\pf
Let $H$ be a subgroup of a hyperbolic group $G$.
By Theorem~\ref{metric}(2), there exist constants $A> 0$ and $B\ge 0$
such that $d(a,b)\leq A\,\hat{d}(a,b)+B$ for all $ a, b\in G$.
Hence $(G, \hat{d})$ is uniformly locally finite and
the $H$-action on $(G,\hat{d})$ is proper.
By Theorem~\ref{metric},  $(G, \hat{d})$ is weakly geodesic  and
strongly bolic,
and the $H$-action on $(G,\hat{d})$ is isometric. By a theorem of
P. de la Harpe and P. Jolissaint, $H$ has property
RD~\cite{h,jp}.
Now Theorem 19 implies
Theorem 20.   
\qed

Theorem 20  has been proved independently 
by Vincent Lafforgue using a different and elegant method~\cite{L4}.

The following result is a direct consequence of Theorem~\ref{bc}.

\begin{ttt} The Kadison-Kaplansky conjecture holds for any torsion free
subgroup $G$ of a hyperbolic group, i.e. there exists no non-trivial projection
in the reduced group $C^{\ast}$-algebra $C_r^{\ast}(G)$.
\end{ttt}

Recall that an element $p$ in $C_r^{\ast}(G)$ is said to be a projection
if $p^{\ast}=p, p^2 =p$. A projection in $C^{\ast}_r(G)$ is said
to be non-trivial if $p\neq 0, 1$. It is well known that the Baum-Connes
conjecture for a torsion free discrete group $G$ implies
the Kadison-Kaplansky conjecture for $G$~\cite{BC, BCH}. 

Michael
Puschnigg has independently proved Theorem 21 using a beautiful
 local cyclic
homology method~\cite{p}. 
 Ronghui Ji has previously proved that there exists no
non-trivial idempotent in the Banach  algebra $\ell^1 (G)$ for any
torsion free hyperbolic group~\cite{J}.

\vspace{0.5cm}

\noindent {\small \sc University of South Alabama\\
Dept of Mathematics and Statistics, ILB 325\\
Mobile, AL 36688, USA}\\
{\tt mineyev@math.usouthal.edu}\\
{\tt http://www.math.usouthal.edu/\~{ }mineyev/math/} 

\vspace{0.5cm}
\noindent {\small \sc Vanderbilt University\\
Department of Mathematics\\
1326 Stevenson Center\\
Nashville, TN 37240}, USA\\
{\tt gyu@math.vanderbilt.edu}
\end{document}